\documentclass{article}
\usepackage{amssymb}
\usepackage{amsthm}
\usepackage{amsmath}
\usepackage[dvips]{epsfig}
\usepackage{graphicx}
\usepackage{color}
\usepackage{url}

\newtheorem{theorem}{\bf Theorem}[section]

\newcommand{\PT}{{\cal PT}}
\newcommand{\T}{{\cal T}}
\newcommand{\p}{{\cal P}}

\begin{document}

\title{Global existence of solutions to coupled ${\cal PT}$-symmetric nonlinear Schr\"odinger equations
}

\author{Dmitry E. Pelinovsky$^1,2$ \and Dmitry A. Zezyulin$^3$ \and Vladimir V. Konotop$^3$\\[2mm]
$^1$ Department of Mathematics, McMaster
  University, Hamilton, Ontario, L8S 4K1, Canada\\[2mm]
$^2$  Department of Applied Mathematics, \\
Nizhny Novgorod State
	Technical University, Nizhny Novgorod, Russia\\[2mm]
$^3$ Centro de F\'isica Te\'{o}rica e Computacional
and Departamento de
              F\'isica,  \\ Faculdade de Ci\^encias, 
              Universidade de Lisboa, \\Avenida Professor
              Gama Pinto 2, Lisboa 1649-003, Portugal
}



\maketitle

\begin{abstract}
We study a system of two coupled  nonlinear Schr\"{o}dinger equations, where
one equation includes gain and the other one includes losses. Strengths of the gain and the loss
are equal, i.e.,   the resulting system is   parity-time ($\PT$) symmetric. The model includes
both linear and nonlinear couplings, such that when all nonlinear coefficients are equal,
the system represents the $\PT$-generalization of the Manakov model. In the one-dimensional case,
we prove the existence of a global solution to the Cauchy problem in energy space $H^1$,
such that the $H^1$-norm of the global solution may grow in time. In the Manakov case,
we show analytically that the $L^2$-norm  of the global solution is bounded for all times
and numerically that the $H^1$-norm is also bounded. In the two-dimensional case, we
obtain a constraint on the $L^2$-norm of the initial data that ensures the existence of a global
solution in the energy space $H^1$.
\end{abstract}

\section{Introduction}
\label{intro}
A system of two coupled nonlinear Schr\"odinger (NLS) equations represents the basic model of propagation of weakly dispersive waves  having vectorial nature, i.e., characterized by two components, below designated by $u$ and $v$. This model is relevant for almost all applications of the nonlinear physics. When the system  is supplied by balanced gain and losses, as well as by the linear coupling, characterized respectively by positive coefficients $\gamma>0$ and $\kappa>0$, it reads
\begin{equation}
\label{Manakov}
\begin{array}{l} i u_t = -  u_{xx} + \kappa v + i \gamma u - (g_{11} |u|^2 + g_{12} |v|^2)u, \\
i v_t =-  v_{xx} + \kappa u - i \gamma v - (g_{12} |u|^2 + g_{22} |v|^2) v, \end{array}
\end{equation}
where all the nonlinear coefficients, i.e., $g_{11}$, $g_{22}$ and $g_{12}$ are real.
Model (\ref{Manakov})  describes wave guiding in structures obeying parity-time  ($\PT$)  symmetry~\cite{PT-general}.
Here the parity symmetry is defined by the mapping $\p (u,v)= (v,u)$.
Time reversal operator $\T$ is defined by the map
$$
\T(u(t),v(t))=(\bar{u}(-t), \bar{v}(-t)).
$$
This means that if $g_{11}=g_{22}$, then  the nonlinearity is also $\PT$ symmetric, i.e.,
for any solution $(u(t), v(t))$ defined in the symmetric interval $[-t_0, t_0]$, there
also exists another solution
$$
u_{\PT}(t) = \bar{v}(-t) \quad \mbox{\rm and} \quad v_{\PT}(t) = \bar{u}(-t).
$$
Hereafter an overbar  stands for the complex conjugation.
Notice that the nonlinear coupling coefficients $g_{12}$ in (\ref{Manakov}) are considered to be
equal in the both equations to ensure that the system is Hamiltonian in the absence of gain and loss with $\gamma = 0$.

Coupled NLS equations (\ref{Manakov}) with $\PT$-symmetric nonlinearity $g_{11}=g_{22}$ and
with $g_{12}=0$, introduced independently in~\cite{DribMal1} for the constant gain-and-loss
coefficient and in~\cite{AKMS} for gain and losses localized along the waveguide, can be
viewed as an extension of the nonlinear dimer model~\cite{coupler_1,coupler_2}.
In the subsequent works~\cite{BKM_dark,BDKM} more general nonlinearities, including ones
with $g_{12}\neq 0$, were thoroughly studied. It was found that the model supports
bright solitons~\cite{DribMal1,BDKM}, dark solitons~\cite{BKM_dark}, Peregrin solitons and rogue waves~\cite{BDKM},
and breathers~\cite{barashenkov1}. For all refereed cases it is important that the nonlinearity
obeys $\PT$ symmetry (for the relevance of $\PT$ symmetry of the nonlinearity, see~\cite{ZK}).

System (\ref{Manakov}) is integrable if $\kappa = \gamma = 0$ and $g_{11} = g_{22} = g_{12}$,
representing the well-known Manakov model~\cite{Manakov} (see also \cite{APT} for review of
these integrable systems). It turns out, that this particular type of (all equal) nonlinearities
is peculiar also for the $\PT$-symmetric model which accounts for gain and losses, as well as
for the linear coupling. More generally, the role of the type of nonlinearity (i.e., non $\PT$-symmetric,
$\PT$-symmetric, and Manakov type) can be crucial for the long time behavior of the solution.
Indication on such possibility stems form the recent results on the stability of trajectories
of the nonlinear $\PT$-symmetric dimer. Namely, the nonlinear dimer
\begin{equation}
\label{dimer}
\begin{array}{l}
i u_t = \kappa v + i \gamma u - |u|^2u,
\\
i v_t =\kappa u - i \gamma v -  |v|^2v,
\end{array}
\end{equation}
has trajectories escaping to infinity (as $t \to \infty$) ~\cite{Barash,Susanto,KPT}, while the ${\cal PT}$-symmetric
Manakov dimer, i.e., the system
\begin{equation}
\label{dimerM}
\begin{array}{l} i u_t = \kappa v + i \gamma u - (|u|^2 + |v|^2)u, \\
i v_t =\kappa u - i \gamma v - ( |u|^2 + |v|^2) v,
\end{array}
\end{equation}
has all trajectories bounded~\cite{KZP}, provided that $\gamma < \kappa$.
Note that the later constraint is related to linear stability of the zero
equilibrium in systems (\ref{Manakov}), (\ref{dimer}), and (\ref{dimerM}).
In the theory of $\PT$-symmetric systems, it is referred to as a condition of
unbroken $\PT$ symmetry~\cite{PT-general}.

In the present paper, we address the existence and boundedness of solutions
to the Cauchy problem for the generalized Manakov system (\ref{Manakov})
with different types of nonlinear coefficients. For the situations where
analytical results are not available, we perform numerical simulations.
We also extend the global existence results to the case of two spatial dimensions.

The paper is organized as follows. Main results in one dimension are presented
in Section 2. They consist of two theorems, the first one about global
well-posedness of the general model (\ref{Manakov}) and
the second one is about bounds on the $L^2$-norm of the solution
in the $\PT$-symmetric Manakov case. The two theorems are then proved in
Sections 3 and 4. Section 5 reports numerical results illustrating that
the $H^1$ norm is also bounded in the $\PT$-symmetric Manakov case, but
it may grow to infinity in the $\PT$-symmetric non-Manakov case.
Finally, extensions of these results to the $\PT$-symmetric Manakov system
in two spatial dimensions are considered in the concluding Section 6.

\section{Main results}

Because system (\ref{Manakov}) is semi-linear, it can be rewritten in an integral form
by using Duhamel's principle. Applying the well-known method based on the contraction mapping
principle (see \cite{LP} for details), one can immediately establish existence
and uniqueness of local solutions of the Cauchy problem for system (\ref{Manakov}) in
Sobolev spaces $H^s(\mathbb{R})$ for any $s > \frac{1}{2}$. In what follows,
we are interested in the behavior of solutions in the energy space $H^1(\mathbb{R})$,
which corresponds to $s = 1$ and appears to be most relevant for physical applications.

Our first result establishes existence of {global} solutions in
the energy space for arbitrary choice of  coefficients in the generalized Manakov system (\ref{Manakov}).

\begin{theorem}
	For any $(u_0,v_0) \in H^1(\mathbb{R}) \times H^1(\mathbb{R})$, there exists a unique global solution
	$(u(t),v(t)) \in C(\mathbb{R},H^1(\mathbb{R}) \times H^1(\mathbb{R}))$ of the generalized Manakov
	system (\ref{Manakov}) such that $(u(0),v(0)) = (u_0,v_0)$.  \label{theorem-1}
\end{theorem}

We remark, that explicit examples of global solutions, including bright solitons,
can be readily obtained  in the case $g_{11} = g_{22}$ and $\gamma < \kappa$, since   the substitution~\cite{DribMal1,BDKM}
$v = e^{ i \delta}u$ with $\delta$ defined by $\sin\delta=-\gamma/\kappa$, reduces
system (\ref{Manakov})  to the integrable nonlinear Schr\"odinger equation
\begin{equation}
i u_t = -u_{xx} + \kappa \cos \delta u - (g_{11} + g_{12}) |u|^2 u.
\end{equation}

The $H^1$-norm of a global solution to the generalized Manakov system (\ref{Manakov})
remains finite but may grow as $t \to \infty$. This growth happens typically in the systems with gains and losses, at least
for sufficiently large initial data, see e.g.~\cite{Barash,Susanto,KPT}. It is then quite surprising that the
Manakov system (\ref{Manakov}) with $g_{11} = g_{12} = g_{22}$ has bounded solutions in the $L^2$-norm
if $\gamma < \kappa$ (i.e. under the condition of the unbroken $\PT$ symmetry). The
next theorem formulates the corresponding result.

\begin{theorem}
	Assume $g_{11} = g_{22} = g_{12}$  and $\gamma < \kappa$.
	For any global solution $(u(t),v(t)) \in C(\mathbb{R},H^1(\mathbb{R}) \times H^1(\mathbb{R}))$ of
	Theorem \ref{theorem-1}, there exists a constant $Q_{max} > 0$ such that
	\begin{equation}
    \label{bound-L2-1}
	\sup_{t \in \mathbb{R}}\left( \| u(t) \|^2_{L^2} + \| v(t) \|^2_{L^2} \right) \leq  Q_{max}.
	\end{equation}
    For $\gamma \geq \kappa$, there exists a global solution of Theorem \ref{theorem-1} such that
    \begin{equation}
       \label{bound-L2-2}
	\lim_{t \to \infty} \left( \| u(t) \|^2_{L^2} + \| v(t) \|^2_{L^2} \right) = \infty.
	\end{equation}
	\label{theorem-2}
\end{theorem}

Note that  Theorem \ref{theorem-2} generalizes the
result of our previous work \cite{KZP} devoted to the integrable ${\cal PT}$-symmetric
dimer (\ref{dimer}) that corresponds to the $x$-independent solutions of system (\ref{Manakov}).
In comparison, the result of Theorem \ref{theorem-1} is novel and the proof is based on
deriving an apriori energy estimates for the $H^1$-norm of the local solutions
and applying Gronwall's inequality.

Also note that the results of Theorems \ref{theorem-1} and \ref{theorem-2} do not clarify
if the $H^1$-norm of the solution is also globally bounded in the case $g_{11} = g_{22} = g_{12}$ and $\kappa < \gamma$.
Numerical simulations reported below in Sec.~\ref{sec:num} support the conjecture that the $H^1$-norm does remain bounded in the dynamics
of the $\PT$-symmetric Manakov system but may grow to infinity in the $\PT$-symmetric non-Manakov case.

\section{Proof of Theorem \ref{theorem-1}}
\label{sec:proof1}

Consider a local solution
\begin{equation}
\label{loc-sol}
(u,v) \in C([-t_0,t_0],H^{1}(\mathbb{R}) \times H^{1}(\mathbb{R}) )
\end{equation}
to the generalized Manakov system (\ref{Manakov}) in the energy space, which exists
for some $t_0 > 0$ due to the standard contraction mapping principle \cite{LP}.
Let us recall the conserved quantities for the Hamiltonian version
of the generalized Manakov equations (\ref{Manakov}) with $\gamma = 0$, which are
well defined in the energy space:
\begin{equation}
\label{power}
Q(t) := \int_{\mathbb{R}} (|u|^2 + |v|^2) dx
\end{equation}
and
\begin{eqnarray}
\nonumber
E(t) & := & \int_{\mathbb{R}} \left( |u_x|^2 + |v_x|^2 +\kappa(\bar{u}v + {u}\bar{v}) \right. \\
\label{energy}
& \phantom{t} & \phantom{text} \left. -
\frac{g_{11}}{2}|u|^4 - \frac{g_{22}}{2}|v|^4 - g_{12} |u|^2 |v|^2 \right) dx.
\end{eqnarray}
If $\gamma \neq 0$, these integral quantities are no longer constant in time. Therefore,
we shall first establish the balance equations for the rate of change of these quantities.
In what follows, we use notations $\| u \|_{L^p}$ for the $L^p$-norm of $u : \mathbb{R} \mapsto \mathbb{C}$.
If $u$ also depends on $t$, we do not usually write the variable $t$ explicitly.

From system (\ref{Manakov}), it is clear that the local solution (\ref{loc-sol})
satisfies $(u_t,v_t) \in C([-t_0,t_0],H^{-1}(\mathbb{R}) \times H^{-1}(\mathbb{R}))$.
Using the duality of Sobolev spaces $H^1(\mathbb{R})$ and $H^{-1}(\mathbb{R})$, we obtain
\begin{equation}
\label{eq:dS0}
\frac{d Q}{dt} = 2\gamma \left( \| u \|_{L^2}^2 - \| v \|_{L^2}^2 \right) \leq 2\gamma Q(t).
\end{equation}
Using Gronwall's inequality for (\ref{eq:dS0}), we deduce that $Q(t)$ does not blow up in a finite time and
grows not faster than an exponential function:
\begin{equation}
\label{bound-l2}
Q(t) \leq Q(0) e^{2\gamma |t|}, \quad t \in [-t_0,t_0].
\end{equation}

The balance equation for the energy $E(t)$ is not defined if we only work
with the local solution (\ref{loc-sol}).
However, we can approximate $H^1$-solutions of system (\ref{Manakov})
by sequences of local solutions of system (\ref{Manakov}) in $H^3(\mathbb{R})$,
for which the energy balance equation is derived using integration
by parts. Then, we take the limits for the sequences and obtain
the following balance equation for the local solutions in $H^1(\mathbb{R})$:
\begin{equation}
\label{energy-balance}
\frac{dE}{dt} = 2 \gamma \int_{\mathbb{R}} \left( |u_x|^2 - |v_x|^2 - g_{11}|u|^4 + g_{22} |v|^4 \right) dx.
\end{equation}
We shall now control the homogeneous $H^1(\mathbb{R})$ norm of the solution by using the integral quantity
\begin{equation}
\label{def-D}
D(t) := \| u_x \|_{L^2}^2 + \| v_x \|_{L^2}^2.
\end{equation}
Note that $D(t)$ is a part of $E(t)$ but it does not satisfy a nice balance equation, for which
the nonlinear terms cancel out. In other words, the energy balance equation for $D(t)$ cannot be estimated
in terms of $D(t)$ only. For this purpose, we use $E(t)$ and, therefore, we shall first establish a correspondence between
$D(t)$ and $E(t)$.

Let us recall the Gagliardo--Nirenberg inequality \cite{LP}.
There exists a positive constant $C_{GN}$ such that for every $f \in H^1(\mathbb{R})$,
we have
\begin{equation}
\label{GNinequality}
\| f \|_{L^4}^4 \leq C_{GN} \| f \|_{L^2}^3 \| f_x \|_{L^2}.
\end{equation}
Using the Gagliardo--Nirenberg inequality (\ref{GNinequality}), we estimate $D(t)$ from above:
\begin{eqnarray}
\nonumber
D(t) & = & E(t) + \int_{\mathbb{R}} \left(\frac{g_{11}}{2}|u|^4 + \frac{g_{22}}{2}|v|^4 + g_{12} |u|^2 |v|^2 -
\kappa(\bar{u}v + {u}\bar{v})\right) dx \\
\nonumber
& \leq & |E(t)| +  g \left(\|u\|_{L^4}^4+\|v\|_{L^4}^4\right)  + \kappa Q(t) \\
\label{eq:D2}
& \leq & |E(t)| + g C_{GN} Q(t)^{3/2} D(t)^{1/2}+\kappa Q(t),
\end{eqnarray}
where $g := \max \{|g_{11}|, |g_{22}|, |g_{12}|\}$. Recall that both $\kappa$ and $\gamma$
are considered to be non-negative for simplicity. From inequality (\ref{eq:D2}), we obtain
\begin{eqnarray}
\label{eq:D}
D(t)^{1/2} \leq \frac{g}{2} C_{GN} Q(t)^{3/2} + \sqrt{|E(t)| + \frac{g^2}{4} C^2_{GN} Q(t)^3+\kappa Q(t)}.
\end{eqnarray}

The energy balance equation (\ref{energy-balance}) can now be estimated again with
the use of the Gagliardo--Nirenberg inequality (\ref{GNinequality}):
\begin{eqnarray}
\nonumber
\left| \frac{dE}{dt} \right| & \leq & 2 \gamma \left( D(t) + g \left(\|u\|_{L^4}^4+\|v\|_{L^4}^4\right) \right) \\
& \leq & 2 \gamma \left(D(t)^2 + g   C_{GN}  Q(t)^{3/2} D(t)^{1/2} \right),\label{eq:E}
\end{eqnarray}
and using (\ref{eq:D2}), (\ref{eq:D}), and the elementary inequality $2 |ab| \leq |a|^2 + |b|^2$ in order to bound $D(t)$,
we obtain
\begin{eqnarray}
\left| \frac{dE}{dt} \right| \leq
4\gamma \left( |E(t)| + \kappa Q(t) + \frac{9}{8} g^2 C^2_{GN} Q(t)^3 \right).\label{eq:E2}
\end{eqnarray}
Using that $E(t) = E(0) + \int_0^t E'(\tau) d\tau$ and integrating (\ref{eq:E2}), we obtain
the integral inequality
\begin{eqnarray}
|E(t)| \leq |E(0)| +  4 \gamma \int_0^t \left( |E(\tau)| + Z(\tau) \right) d \tau,\label{eq:E3}
\end{eqnarray}
where
\begin{eqnarray}
Z(t) := \kappa Q(t) + \frac{9}{8} g^2 C^2_{GN} Q(t)^3.
\end{eqnarray}
  Notice that because of the bound (\ref{bound-l2}), function $Z(t)$ is defined for all $t \in [-t_0,t_0]$ and
does not grow faster than an exponential function. Using Gronwall's inequality for (\ref{eq:E3}),
we obtain
\begin{equation}
\label{eq-E4}
|E(t)| \leq  \left( |E(0)|  +  4 \gamma \int_0^t Z(\tau) d\tau \right) e^{4\gamma |t|}, \quad t \in [-t_0,t_0].
\end{equation}
Therefore $|E(t)|$ does not blow up in a finite time and grows not faster than an exponential function.
The same is true  for $D(t)$ thanks to the bound (\ref{eq:D}).

Extending the local solution of the generalized Manakov system (\ref{Manakov}) in $H^1(\mathbb{R})$
to larger values of $t_0$ by using a priori energy inequalities (\ref{bound-l2}), (\ref{eq:D}), and (\ref{eq-E4})
and applying the standard continuation arguments, we conclude the proof of Theorem \ref{theorem-1} on the existence of
the unique global solution of the generalized Manakov system (\ref{Manakov}) in $H^1(\mathbb{R})$.

\section{Proof of Theorem \ref{theorem-2}}
\label{sec:proof2}
For the global solution of Theorem \ref{theorem-1}, we can introduce the integral Stokes  variables
\begin{eqnarray*}
	S_1(t)  &:=&  \int_{\mathbb{R}} (\bar{u} v + u \bar{v}) dx, \\
	S_2(t)  &:=&  i \int_{\mathbb{R}} (\bar{u} v - u \bar{v}) dx,\\
	S_3(t)  &:=&  \int_{\mathbb{R}} (|u|^2 - |v|^2) dx.
\end{eqnarray*}
Let us consider the case $g_{11} = g_{22} = g_{12}$.
We compute temporal derivatives of the integral Stokes variables and obtain
\begin{eqnarray}
\label{Stokes-system}
\frac{dS_1(t)}{dt} = 0, \quad \frac{dS_2(t)}{dt} = 2\kappa S_3, \quad \frac{dS_3(t)}{dt} = 2\gamma Q - 2\kappa S_2,
\end{eqnarray}
where $Q$ is defined by (\ref{power}) and satisfies the balance equation (\ref{eq:dS0}).
Note that all nonlinear terms cancel out in the system of evolution equations (\ref{Stokes-system})
for the integral Stokes variables.

It follows from system (\ref{Stokes-system}) that $S_1$ and $C :=\kappa Q-\gamma S_2$ are constants of motion.
Excluding $S_2$ and $S_3$ from system (\ref{Stokes-system}), we obtain the linear oscillator equation for $Q(t)$
\begin{eqnarray}
\label{lin-osc}
\frac{d^2 Q}{d t^2} + 4\omega^2  Q = 4 \kappa C = {\rm const}.
\end{eqnarray}
where we have introduced
\begin{eqnarray}
\label{omega}
\omega=\sqrt{\kappa^2 - \gamma^2}
\end{eqnarray}
The linear oscillator equation (\ref{lin-osc}) immediately implies that if
$\gamma < \kappa$, then $Q(t)$ oscillates with the frequency $2 \omega$
and therefore, can be globally bounded
for all $t \in \mathbb{R}$ by a constant $Q_{max}>0$.
Indeed, from the linear equation (\ref{lin-osc})
with the initial conditions $Q_0 = Q(0)$ and $P_0 = \dot{Q}(0) = 2\gamma S_3(0)$,
we compute
\begin{equation}
Q(t) = \frac{\kappa C}{\omega^2} + A_1 \cos(2\omega t) + A_2 \sin(2\omega t),
\end{equation}
where integration constants are given by
\begin{equation}
A_1 = Q_0 - \frac{\kappa^2 C}{\omega^2}, \quad \mbox{and}\quad A_2 = \frac{P_0}{2 \omega}.
\end{equation}
Hence  $Q(t) \leq Q_{max}$, where the sharp upper bound  $Q_{max}$ is computed from the initial conditions by
\begin{equation}
\label{Qmax}
 Q_{max} = \frac{\kappa C}{\omega^2} + \sqrt{A_1^2 + A_2^2}
\end{equation}
This argument proves bound (\ref{bound-L2-1}) of Theorem \ref{theorem-2}.

If $\gamma \geq \kappa$, there exists solutions to Eq. (\ref{lin-osc}) with $\omega^2 \leq 0$
such that $Q(t)$ grows unboundedly as $t \to \infty$. Therefore, the $L^2$-norm of global solutions
of the $\PT$-symmetric Manakov system with  $\gamma\geq \kappa$ grows to infinity,
according to bound (\ref{bound-L2-2}). This concludes the proof of Theorem \ref{theorem-2}.

Note that the linear oscillator equation (\ref{lin-osc}) is derived in our previous work \cite{KZP}
for the $x$-independent solutions by using the local Stokes variables. Furthermore,
we note here that for $x$-independent solutions, the ${\cal PT}$-symmetric
Manakov dimer (\ref{dimerM}) can be reduced to the linear system. Indeed, considering
$$
u(t) = \tilde{u}(t) e^{i \varphi(t)}, \quad v(t) = \tilde{v}(t) e^{i \varphi(t)},
$$
and setting $\varphi(t) = \int_0^t \left( |\tilde{u}|^2 + |\tilde{v}|^2 \right) dt$, we obtain the system of linear equations
\begin{equation}
\label{dimer-linear}
\begin{array}{l}
i \tilde{u}_t = \kappa \tilde{v} + i \gamma \tilde{u}, \\
i \tilde{v}_t =\kappa \tilde{u} - i \gamma \tilde{v}, \end{array}
\end{equation}
which is stable for $\gamma < \kappa$ and unstable for $\gamma \geq \kappa$. Although these results
are rather trivial for $x$-independent solutions to the ${\cal PT}$-symmetric Manakov dimer (\ref{dimerM}),
we have shown here that these results can be extended to bound the $L^2$-norm of $x$-dependent solutions
to the ${\cal PT}$-symmetric Manakov system (\ref{Manakov}) with $g_{11} = g_{22} = g_{12}$ in the energy space.

\section{Numerical simulations}
\label{sec:num}
The results of Theorems \ref{theorem-1} and \ref{theorem-2} do not clarify if the $H^1$-norm
of global solutions   of the $\PT$-symmetric  system (\ref{Manakov}) remains
bounded for $\gamma < \kappa$. In order to obtain additional insight into this issue,
we have undertaken a set of numerical simulations.
Our numerical results indicate that the $H^1$-norm of the global solutions
does remain bounded for the Manakov system  ($g_{11} = g_{22} = g_{12}$) with
unbroken $\PT$ symmetry ($\gamma < \kappa$). On the other hand, in the $\PT$-symmetric
non-Manakov case with $g_{11} = g_{22} \neq g_{12}$, unbounded solutions  exist
even if $\gamma < \kappa$. The analytical proof of these conjectures remains opened for further studies.

More specifically, we run a set of numerical simulations of the dynamics governed by
the generalized Manakov system (\ref{Manakov}) which is solved  in the finite domain
$[-l, l]$ with zero boundary conditions, $u(\pm l,t)=v(\pm l,t)=0$, using a semi-implicit difference scheme.
The numerically obtained solution is used to evaluate quantities $D(t)$ and  $Q(t)$ which
characterize the $H^1$-norm of the global solution. For the initial conditions at $t = 0$,
we use Gaussian beams
\begin{equation}
\label{IC}
u_0(x) = \frac{A}{\pi^{1/4} a^{1/2}} \exp\left\{ -\frac{x^2}{2a^2}\right \}, \quad
v_0(x) = \frac{B}{\pi^{1/4} b^{1/2}} \exp\left\{ -\frac{x^2}{2b^2}\right \}
\end{equation}
truncated on the interval $(-l,l)$, where $A$, $B$, $a$ and $b$ are real constants.
Notice that the  functions $u_0$ and $v_0$ are normalized such that
$$
\int_{-\infty}^\infty |u_0(x)|^2 dx = A^2 \quad \mbox{\rm and} \quad
\int_{-\infty}^\infty |v_0(x)|^2 dx = B^2.
$$
The width $l$ of the computational domain is sufficiently large such that
the numerical value of the  initial conditions  (\ref{IC}) nearly vanish at $x=\pm l$.

\begin{figure} 
 \centering
	\includegraphics[width=0.8\textwidth]{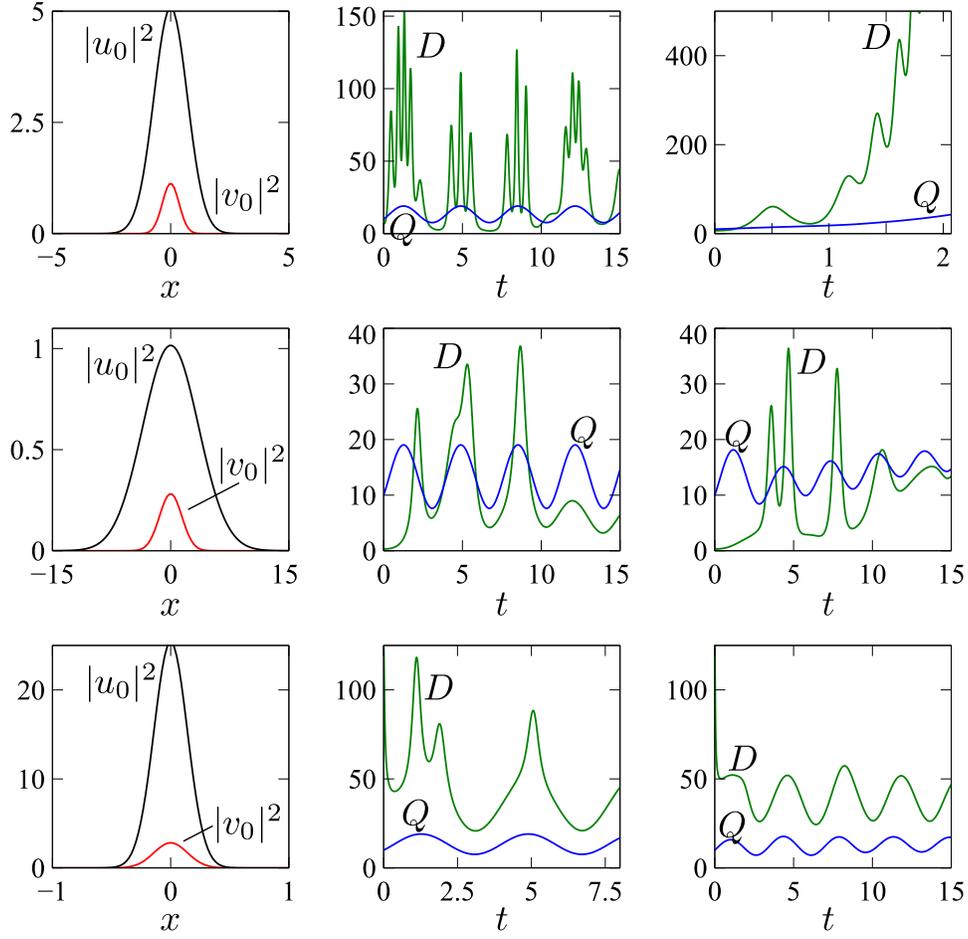}
	 \caption{Results of numerical solution of the Cauchy problem for $\PT$-symmetric   system (\ref{Manakov}). Three rows present results obtained for different initial conditions (\ref{IC}): upper, middle and the lower rows correspond to $(A, a, B,b)$ as follows:  (3, 1, 1, 0.5),  (3, 5, 1, 2), and  (3, 0.2, 1, 0.2). Panels from the left  column show profiles of the initial conditions $|u_0(x)|^2$ (black lines) and $|v_0(x)|^2$ (red lines). Panels from the central column show  dependencies $D(t)$ (green lines, labeled as $D$) and   $Q(t)$ (blue  lines, labeled as $Q$) found for the Manakov case $g_{11}=g_{22} = g_{12}=1$. Panels from the right column show the same data, obtained from the same initial conditions as in the middle column,  but for the  {non-Manakov} case   $g_{11}=g_{22} =1$ and $g_{12}=1/2$. For all panels, $\kappa=1$ and $\gamma=1/2$, i.e., $\PT$ symmetry is unbroken.
}
\label{fig_conj}
\end{figure}

We solved   the Cauchy problem in the finite domain $[-l, l]$ for   several different sets of initials conditions $(A, a, B,b)$, as well as for different combinations of $\gamma$ and $\kappa$ (with $\gamma < \kappa$) and for different $l$. In Fig.~\ref{fig_conj}  we present results of   numerical runs for three different initial conditions.
Figure~\ref{fig_conj}  addresses the case of the Manakov nonlinearity with $g_{11} = g_{22} = g_{12} = 1$
(panels in the middle column) and the non-Manakov case with $g_{11} = g_{22} = 1$ and $g_{12} = 1/2$
(panels in the right column).

For all sets of the initial conditions (\ref{IC}) and for all $l$ we have tested,
we  observed that in the case $g_{11}=g_{22}=g_{12}$, the numerically obtained $H^1$ norm
remains bounded (see panels in the middle columns of Fig.~\ref{fig_conj}). This allows us
to conjecture on boundness of the $H^1$-norm of all global solutions to the $\PT$-symmetric
Manakov system (\ref{Manakov}) with  $g_{11}=g_{22}=g_{12}$ and $\gamma < \kappa$. On the other hand, unbounded numerical
solutions exist in the non-Manakov case (see the upper panel in the right column of Fig.~\ref{fig_conj}).   In all panels showing dynamics of the Manakov system, we observe the oscillatory
behavior of $Q(t)$ which is characterized by the period (see (\ref{lin-osc}) and (\ref{omega})),
\begin{equation}
\frac{\pi}{\omega} = \frac{\pi}{\sqrt{\kappa^2-\gamma^2}}\approx 3.6276.
\end{equation}

Since  our main interest is  related to the existence (or non-existence)
of globally bounded   solutions, in Fig.~\ref{fig_conj}  we show the results obtained for a situation
when the input beam in the equation with gain  has the intensity much larger than the intensity of
the input beam for the  NLS equation with dissipation (c.f. black and red curves in the left columns
of Fig.~\ref{fig_conj}). The opposite situation when the intensity of the input beam for the
lossy equation is  much larger than that for the equation with gain was also tested in our
numerics with the same  conclusions on boundedness and unboundedness of solutions.

\section{Extension to the two-dimensional Manakov system}

In this concluding section, we generalize the proof of the global existence in the energy space for the Manakov system, i.e. for (\ref{Manakov}) with all equal nonlinear coefficients, presented above, to the case of two spatial dimensions, $x\in\mathbb{R}^2$. In this context we note that the analysis of  global existence and blow-up of solutions in a Hamiltonian counterpart of this system, it corresponds to $\gamma=0$, was recently presented in \cite{JW13}.

Respectively we shall rewrite the system (\ref{Manakov}) in two spatial dimensions
with the Manakov-type nonlinearity:
\begin{equation}
\label{Manakov2D}
\begin{array}{l}
i u_t = - \Delta u  + \kappa v + i \gamma u - (|u|^2 +  |v|^2)u,
\\
i v_t =-  \Delta v + \kappa u - i \gamma v - (|u|^2 + |v|^2) v.
\end{array}
\end{equation}
Here $\Delta$ is the two-dimensional Laplacian and all the nonlinear coefficients
$g_{11}$, $g_{22}$ and $g_{12}$ are set to unity.


We shall now formulate and prove the global existence result for the $\PT$-symmetric Manakov system
(\ref{Manakov2D}) by generalizing the ideas of Theorems~\ref{theorem-1} and \ref{theorem-2}.
This generalization   is possible due to the fact that the $L^2$-norm of a global  solution
to system (\ref{Manakov2D})  is bounded for $\gamma < \kappa$
by the constant $Q_{max}$, which is determined in terms of the initial conditions
in $H^1(\mathbb{R}^2)$ by the explicit expression (\ref{Qmax}).
This global bound allows us to obtain an a priori bound for the energy
\begin{equation*}
E(t) := \int_{\mathbb{R}^2} \left( |\nabla u|^2 + |\nabla v|^2 +\kappa(\bar{u}v + {u}\bar{v})-
\frac{1}{2}|u|^4 - \frac{1}{2}|v|^4 - |u|^2 |v|^2 \right) dx.
\end{equation*}
Note here that the local solution in the energy space
\begin{equation}
\label{loc-sol-2}
(u(t),v(t)) \in C([-t_0,t_0],H^1(\mathbb{R}^2) \times H^1(\mathbb{R}^2)),
\end{equation}
is not derived by the standard contraction mapping principle for the semi-linear equations.
Nevertheless, a modified contraction method can be developed with the use of Stritcharz
estimates as in \cite{LP}. Therefore, we use here the local solution (\ref{loc-sol-2})
obtained from this modified contraction method.

Next, let us recall the Gagliardo-Nirenberg inequality in the space of two dimensions \cite{LP}.
There exists a positive constant $C_{GN}$ such that for every $f\in H^1(\mathbb{R}^2)$, we have
\begin{equation}
\label{GN2D}
\| f\|_{L^4}^4 \leq C_{GN} \|\nabla f \|_{L^2}^2    \|f\|_{L^2}^2.
\end{equation}
Moreover, it was shown by Weinstein \cite{W83} that $C_{GN}$ can be found
sharply in terms of the ground state solution $R$ of the stationary equation
$$
\Delta R - R + R^3 = 0,
$$
namely, $C_{GN} = 2/\|R\|_{L^2} \approx 0.171$.
In physics,  the solution $R(x)$ is also known as the Townes soliton~\cite{T64}.
Now we can formulate a global existence result for the two-dimensional
$\PT$-symmetric Manakov system (\ref{Manakov2D}).

\begin{theorem}
    Consider  the  two-dimensional $\PT$-symmetric Manakov system (\ref{Manakov2D}) with $\gamma < \kappa$.
    If the initial conditions $(u_0, v_0) \in H^1(\mathbb{R}^2) \times H^1(\mathbb{R}^2)$ satisfy the constraint
    \begin{equation}
    \label{condition2D}
    Q_{max} < \frac{1}{2}  \|R\|_{L^2}^2
    \end{equation}
    where $Q_{max}$  is defined by (\ref{Qmax}), then there exists a unique global solution
    $(u(t),v(t)) \in C(\mathbb{R},H^1(\mathbb{R}^2) \times H^1(\mathbb{R}^2))$ such that $u(0) = u_0$ and  $v(0) = v_0$.
    On the other hand, if $Q_{max} \geq \|R \|_{L^2}^2$, there exists initial conditions
    $(u_0, v_0) \in H^1(\mathbb{R}^2) \times H^1(\mathbb{R}^2)$, such that local solutions
    (\ref{loc-sol-2}) blow up in a finite time, that is, there exists $T_0 > 0$ such that
    \begin{equation}
    \label{blowup2D}
    \lim_{t \uparrow T} \left( \| u \|_{H^1}^2 + \| v \|_{H^1}^2 \right) = \infty.
    \end{equation}
	\label{theorem-2D}
\end{theorem}

We shall now sketch the proof  of Theorem~\ref{theorem-2D} by following the ideas of the proof of Theorem~\ref{theorem-1}.
Defining $D(t)$ as in (\ref{def-D}), we employ   the Gagliardo-Nirenberg inequality (\ref{GN2D}) and estimate $D(t)$ from above as
\begin{eqnarray*}
D(t) & \leq & |E(t)| + \| u \|_{L^4}^4 +   \| v \|_{L^4}^4 + \kappa Q(t) \\
& \leq & |E(t)| + C_{GN} D(t) Q(t) + \kappa Q(t).
\end{eqnarray*}
For  all $t\in [-t_0, t_0]$ we have $Q(t) \leq Q_{max}$. If the condition (\ref{condition2D}) is met (recall that
$C_{GN} = 2/\| R \|_{L^2}^2$), then $D(t)$ is controlled by $|E(t)|$ from above as
\begin{equation}
D(t) \leq \frac{|E(t)| + \kappa Q_{max}}{1 - C_{GN} Q_{max}}.
\end{equation}
Using then the balance equation for $dE/dt$ as in (\ref{energy-balance}), we obtain
for all $t\in [-t_0, t_0]$
\begin{eqnarray*}
\left| \frac{dE }{dt}\right| & \leq & 2 \gamma \left( D(t) + \| u \|_{L^4}^2 + \| v \|^4_{L^4} \right) \\
& \leq & 2 \gamma \left( D(t) + C_{GN} D(t) Q(t) \right) \\
& \leq & 2 \gamma \frac{1 + C_{GN} Q_{max}}{1 - C_{GN} Q_{max}} \left( |E(t)| + \kappa Q_{max} \right).
\end{eqnarray*}
Integrating the latter inequality similarly to what is done in the derivation of (\ref{eq:E3}) and (\ref{eq-E4}), we obtain
a global (exponentially growing) bound on the $H^1$-norm of the local solution.  Hence, the local solution (\ref{loc-sol-2})
can be extended to larger values of $t_0$, and finally, globally for all $t \in \mathbb{R}$.

To show the other part of Theorem \ref{theorem-2D} concerning the finite-time blow-up (\ref{blowup2D}),
we note the following trick, which is explored in the recent works \cite{DribMal1,BDKM,DFKZ}.
Using the substitution $v = e^{i \delta}u$ with $\delta$ defined by $\sin \delta=-\gamma/\kappa$,
we reduce the $\PT$-symmetric Manakov system (\ref{Manakov2D}) with $\gamma < \kappa$ to a scalar NLS equation
\begin{eqnarray}
\label{NLS2D}
iu_t=- \Delta u-2|u|^2u + \omega  u.
\end{eqnarray}
Solutions of equation (\ref{NLS2D}) conserve the $L^2$-norm and may blow up in a finite time
if $\|u_0\|_{L^2}^2 \geq\|R\|_{L^2}^2/2$~\cite{W83}. The corresponding blowing-up solution
of the $\PT$-symmetric Manakov system (\ref{Manakov2D}) also conserves the $L^2$ norm, which then
satisfies $Q(t) = Q_{max} = 2 \| u_0 \|_{L^2}^2 \geq \|R\|_{L^2}^2$. This concludes the proof of Theorem \ref{theorem-2D}.

We note that Theorem \ref{theorem-2D} does not clarify if solutions to the $\PT$-symmetric
Manakov system (\ref{Manakov2D}) may blow up in a finite time if the initial data corresponds
to $Q_{max}$ fitting in the interval
$$
\frac 12 \|R\|_{L^2}^2 \leq Q_{max}<\|R\|_{L^2}^2.
$$
This question remains open for further studies.

\section{Acknowledgements}
The authors thank Vadim Vekslerchik for helpful discussion at an earlier stage of the project. The work of D.P. is supported by the Ministry of Education and Science of Russian Federation (the base part of the state task No. 2014/133).
DAZ and VVK acknowledge support of FCT (Portugal) under the grants  PEst-OE/FIS/UI0618/2014 and  PTDC/FIS-OPT/1918/2012.



\end{document}